\documentclass[11pt,dvips]{article}
\usepackage{theorem}
\usepackage{epsfig}
\usepackage{amsfonts}
\usepackage{amssymb}
\usepackage{latexsym}

\newcommand{\spc}{\vspace{2mm}}

\newcommand{\Z}{\mathbb{Z}}
\newcommand{\F}{\mathbb{F}}

\usepackage{url,amsmath}

\begin{document}

\title{GROUP THEORY IN CRYPTOGRAPHY}
\author{SIMON R. BLACKBURN, CARLOS CID and CIARAN MULLAN\\
Department of Mathematics, Royal Holloway, University of
London\\
Egham, Surrey TW20 0EX, United Kingdom\\
\texttt{\{s.blackburn,carlos.cid,c.mullan\}@rhul.ac.uk}}
\maketitle

\begin{abstract}
This paper is a guide for the pure mathematician who would like to
know more about cryptography based on group theory. The paper gives a
brief overview of the subject, and provides pointers to good
textbooks, key research papers and recent survey papers in the area.
\end{abstract}

\section{Introduction}
\label{sec:introduction}

In the last few years, many papers have proposed cryptosystems based
on group theoretic concepts. Notes from a recent advanced course on
the subject by Myasnikov, Shpilrain and Ushakov have recently been
published as a monograph~\cite{Myasnikov}, and a textbook (with a
rather different focus) by Gonz\'alez Vasco, Magliveras and
Steinwandt~\cite{GonzalezVasco} is promised in 2010. Group-based
cryptosystems have not yet led to practical schemes to rival RSA and
Diffie--Hellman, but the ideas are interesting and the different
perspective leads to some worthwhile group theory. The cryptographic
literature is vast and diverse, and it is difficult for a newcomer to
the area to find the right sources to learn from. (For example, there
are many introductory textbooks aimed at the mathematical audience
that introduce RSA. How many of these textbooks hint that the basic
RSA scheme is insecure if refinements such as message padding are not
used? For a discussion of these issues, see Smart~\cite[Chapters 17,18
and~20]{Smart}, for example.) Our paper will provide some pointers to
some sources that, in our opinion, provide a good preparation for
reading the literature on group-based cryptography; the paper will
also provide a high level overview of the subject. We are assuming
that our reader already has a good knowledge of group theory, and a
passing acquaintance with cryptography: the RSA and Diffie--Hellman
schemes have been met before, and the difference between a public key
and a symmetric key cipher is known.

The remainder of the paper is structured as follows. In
Section~\ref{sec:cryptography} we review some of the basic concepts of
cryptography we will need. In
Section~\ref{sec:cryptographyusinggroups} we introduce some of the
most widely studied schemes in group-based cryptography, and in
Section~\ref{sec:cryptanalysis} we sketch attacks on these schemes. In
all these sections, we cite references that provide more
details. Finally, in Section~\ref{sec:nextsteps}, we touch on some related
areas and give suggestions as to where to search for current papers
and preprints in the subject.

\section{Cryptography Basics}
\label{sec:cryptography}

There are innumerable books on cryptography that are written for a
popular audience: they almost always take a historical approach to the
subject. For those looking for a definitive historical reference book,
we would recommend Kahn~\cite{Kahn} for an encyclopedic and
beautifully written account.

Technical introductions to the area written for a mathematical
audience tend to concentrate (understandably, but regrettably from the
perspective of a cryptographer) on the areas of cryptography that have
the most mathematical content. Stinson~\cite{Stinson} is a well-written 
introduction that avoids this pitfall. Another good reference
is Smart~\cite{Smart}, which has the advantage of being available
online for free. Once these basics are known, we suggest reading a
book that looks at cryptography from the perspective of theoretical
computer science and complexity theory: Katz and
Lindell~\cite{KatzLindell} is a book we very much enjoy. The
theoretical computer science approach has had a major influence on the
field, but is not without its controversial aspects: see
Koblitz~\cite{Koblitz} and responses by Goldreich and
others~\cite{Goldreich}. For readers who insist on falling into the
mathematical pit mentioned above, the book by
Washington~\cite{Washington} on cryptography using elliptic curves is
an excellent follow-up read; elliptic curve based cryptography is
becoming the norm for the current generation of public key
cryptosystems. As we are writing for a mathematical audience, we also
consciously aim to fall into this pit.

A standard model for a cryptographic scheme is phrased as two parties,
Alice and Bob, who wish to communicate securely over an insecure
channel (such as a wireless link, or a conventional phone line). If
Alice and Bob possess information in common that only they know (a
shared secret key) they can use this, together with a symmetric key cipher
such as AES (the Advanced Encryption Standard), to communicate. If
Alice and Bob do not possess a secret key, they execute a protocol
such as the Diffie--Hellman key agreement protocol to create one, or use a
public key cryptosystem such as RSA or ElGamal that does not need a
secret key. Many of the schemes we discuss are related to the
Diffie--Hellman protocol, so we give a brief description of this
protocol as a reminder to the reader.

\spc\textbf{Diffie--Hellman Key Agreement Protocol~\cite{diffie-hellman}.}
Let $G$ be a cyclic group, and $g$ a generator of $G$, where both $g$
and its order $d$ are publicly known. If Alice and Bob wish to create
a shared key, they can proceed as follows:
\begin{enumerate}
\item Alice selects uniformly at random an integer $a \in [2, d-1]$, computes $g^a$, and sends it to Bob.
\item Bob selects uniformly at random an integer $b \in [2, d-1]$, computes $g^b$, and sends it to Alice.
\item Alice computes $k_a=(g^b)^a$, while Bob computes $k_b=(g^a)^b$.
\item The shared key is thus $k = k_a = k_b \in G$.
\end{enumerate}
\spc

The security of the scheme relies on the assumption that, knowing $g
\in G$ and having observed both $g^a$ and $g^b$, it is computationally
infeasible for an adversary to obtain the shared key.  This is known
as the \textbf{Diffie--Hellman Problem (DHP)}. The Diffie--Hellman
problem is related to a better known problem, the Discrete Logarithm
Problem:

\spc \textbf{Discrete Logarithm Problem (DLP).}
Let $G$ be a cyclic group, and $g$ a generator of $G$. Given $h \in
G$, find an integer $t$ such that $g^t = h$.
\spc

Clearly, if the DLP is easy then so is the DHP and thus the
Diffie--Hellman key agreement protocol is insecure. So, as a minimum
requirement, we are interested in finding difficult instances of the
DLP. It is clear that difficulty of the DLP depends heavily on the way
the group $G$ is represented, not just on the isomorphism class of
$G$. For example, the DLP is trivial if $G = \Z /d\Z$ is the additive
group generated by $g=1$. However, if $G$ is an appropriately chosen
group of large size, the DLP is considered computationally
infeasible. In practice, one often uses $G =
\F^*_{p^l}$ (for appropriately selected prime $p$ and exponent~$l$),
or the group of points of a properly chosen elliptic curve over a
finite field.

Turning from the Diffie--Hellman scheme to the more general model,
there are two points we would like to emphasise:
\begin{itemize}
\item \textbf{Alice and Bob are computers.}
So our aim is to create a protocol that is well-specified enough to be
implemented. In particular, a well specified scheme must describe how
group elements are stored and manipulated; the scheme's description
must include an algorithm to generate any system-wide parameters; it
must be clear how any random choices are made. (This last point is
especially critical if we are choosing elements from an infinite set,
such as a free group!) Moreover, the protocol should be efficient; the
computational time required to execute the protocol is critical, but
so are: the number of bits that need to be exchanged between Alice and
Bob; the number of passes (exchanges of information) that are needed
in the protocol; the sizes of keys; the sizes of system parameters.
\item \textbf{Security is a very subtle notion.}
For the last 100 years, it has become standard for cryptographers to
assume that any eavesdropper knows everything about the system that is
being used apart from secret keys and the random choices made by
individual parties. (Claude Shannon~\cite[Page~662]{Shannon} phrased this as
`The enemy knows the system being used'; the phrase `The enemy knows
the system' is known as Shannon's maxim.). But modern security is
often much more demanding. For example, in the commonly studied
IND-CCA2 model, we require that an eavesdropper cannot feasibly guess
(with success probability significantly greater than $0.5$) which of
two messages has been encrypted, when they are presented with a single
challenge ciphertext that is an encryption of one of the
messages. This should even be true when the eavesdropper can choose the
two messages, and is allowed to request the decryption of
any ciphertext not equal to the challenge ciphertext. Note that
cryptographers are usually interested in the complexity in the generic
case (in other words, what happens most of the time). Worst case
security estimates might not be useful in practice, as the worst case
might be very rare; even average case estimates might be unduly
distorted by rare but complicated events. See Myasnikov et
al.~\cite{Myasnikov} for a convincing argument on this point in the
context of group-based cryptography.
\end{itemize}

We end the section by making the point that modern cryptography is
much broader than the traditional two party communication model we
have discussed here: there is a thriving community developing the
theory of multi-party communication, using such beautiful concepts as
zero knowledge. See Stinson~\cite[Chapter~13]{Stinson} for an
introduction to zero knowledge, and see the links from Helger Lipmaa's
page~\cite{Lipmaa} for some of the important papers on multi-party
computation.

\section{Cryptography Using Groups}
\label{sec:cryptographyusinggroups}

This section will discuss several ways in which group theory can be
used to construct variants of the Diffie--Hellman key agreement
protocol. Since the protocol uses a cyclic subgroup of a finite group
$G$, one approach is to search for examples of groups that can be
efficiently represented and manipulated, and that possess cyclic
subgroups with a DLP that seems hard. Various authors have suggested
using a cyclic subgroup of a matrix group in this context, but some
basic linear algebra shows that this approach is not very useful: the
DLP is no harder than the case when $G$ is the multiplicative group of
a finite field; see Menezes and Vanstone~\cite{MenezesVanstone} for
more details. Biggs~\cite{Biggs} has proposed representing an abelian
group as a critical group of a finite graph; but
Blackburn~\cite{Blackburn} has shown that this proposal is
insecure. An approach (from number theory rather than
group theory) that has had more success is to consider the group of
points on an elliptic curve, or Jacobians of hyperelliptic curves. See
Galbraith and Menezes~\cite{GalbraithMenezes} for a survey of this
area.

All the proposals discussed above use representations of abelian
(indeed, cyclic) groups. What about non-abelian groups? The first
proposal to use non-abelian groups that we are aware of is due to
Wagner and Magyarik~\cite{WagnerMagyarik} in 1985. (See Gonz\'alez
Vasco and Steinwandt~\cite{GonzalezVascoSteinwandt} for an attack on
this proposal; see Levy-dit-Vehel and
Perret~\cite{LevyditVehelPerret,LevyditVehelPerret2} for more recent
related work.) But interest in the field increased with two high-profile 
proposals approximately ten years ago. We now describe these proposals.

\subsection{Conjugacy and exponentiation}

Let $G$ be a non-abelian group. For $g,x\in G$ we write
$g^x$ for $x^{-1}gx$, the conjugate of $g$ by $x$. The notation suggests
that conjugation might be used instead of exponentiation in
cryptographic contexts. So we can define an analogue
to the discrete logarithm problem:

\spc
\textbf{Conjugacy Search Problem.}
Let $G$ be a non-abelian group. Let $g,h \in G$ be such that $h =
g^x$ for some $x \in G$. Given the elements $g$ and $h$,
find an element $y \in G$ such that $h = g^y$.
\spc

Assuming that we can find a group where the conjugacy search problem
is hard (and assuming the elements of this group are easy to store and
manipulate), one can define cryptosystems that are analogues of
cryptosystems based on the discrete logarithm problem.  Ko et
al.\ proposed the following analogue of the
Diffie--Hellman key agreement protocol.

\spc
\textbf{Ko--Lee--Cheon--Han--Kang--Park Key Agreement Protocol~\cite{ko-et-al}.}
Let $G$ be a non-abelian group, and let $g$ be a publicly known element of
$G$.  Let $A,B$ be commuting subgroups of $G$, so $[a,b] = 1$ for all $a
\in A$, $b \in B$.  If Alice and Bob wish to create a common secret key, they
can proceed as follows:
\begin{enumerate}
\item Alice selects at random an element $a \in A$, computes $g^a = a^{-1} g a$, and sends it to Bob.
\item Bob selects at random an element $b \in B$, computes $g^b = b^{-1} g b$, and sends it to Alice.
\item Alice computes $k_a = (g^b)^a$, while Bob computes $k_b = (g^a)^b$.
\item Since $ab=ba$, we have $k_a = k_b$, as group elements (though their
representations might be different). For many groups, we can use $k_a$
and $k_b$ to compute a secret key. For example, if $G$ has an
efficient algorithm to compute a normal form for a group element, the
secret key $k$ could be the normal form of $k_a$ and $k_b$.
\end{enumerate}
\spc

The interest in the paper of Ko et al.~\cite{ko-et-al} centred on
their proposal for a concrete candidate for $G$ and the subgroups $A$
and $B$, as follows. We take $G$ to be the braid group $B_n$ on $n$
strings (see Artin~\cite{Artin}, for example) which has presentation
\[
B_n=\left\langle \sigma_1,\sigma_2,\ldots ,\sigma_{n-1}\Big\vert\begin{array}{cl}
\sigma_i\sigma_j\sigma_i=\sigma_j\sigma_i\sigma_j&\text{ for }|i-j|=1\\
\sigma_i\sigma_j=\sigma_j\sigma_i&\text{ for }|i-j|\geq 2\end{array}\right\rangle.
\]
Let $l$ and $r$ be integers such that $l+r=n$. Then we take
\begin{align*}
A&=\left\langle \sigma_1,\sigma_2,\ldots,\sigma_{l-1}\right\rangle\text{ and}\\
B&=\left\langle \sigma_{l+1},\sigma_{l+2},\ldots,\sigma_{l+r-1}\right\rangle.
\end{align*}

The braid group is an attractive choice for the underlying group (a
so-called `platform group') in the Ko et al.\ key agreement protocol:
there is an efficient normal form for an element; group multiplication
and inversion can be carried out efficiently; the conjugacy problem
looks hard for braid groups. Note that we have not specified the
cryptosystem precisely. Of course, we have not chosen the values of
$n$, $l$ and $r$. But we have also not specified how to choose the element
$g\in G$ (it emerges that this choice is critical). Finally, since the
subgroups $A$ and $B$ are infinite, it is not obvious how the elements
$a\in A$ and $b\in B$ should be chosen.

\subsection{Computing a common commutator}

The following beautiful key agreement protocol, due to Anshel, Anshel
and Goldfeld~\cite{Anshel}, has an advantage over the Ko et
al.\ protocol: commuting subgroups $A$ and $B$ are not needed.

\spc\textbf{Anshel--Anshel--Goldfeld Key Agreement Protocol~\cite{Anshel}.}
Let $G$ be a non-abelian group, and let elements $a_1,\ldots
,a_k,b_1,\ldots ,b_m\in G$ be public. 
\begin{enumerate}
\item Alice picks a private word $x$ in $a_1,\ldots ,a_k$ and sends $b^x_1,\ldots ,b^x_m$ to Bob.
\item Bob picks a private word $y$ in $b_1,\ldots ,b_m$ and sends $a^y_1,\ldots ,a^y_k$ to Alice.
\item Alice computes $x^y$ and Bob computes $y^x$.
\item The secret key is $[x,y]=x^{-1}y^{-1}xy$.
\end{enumerate}
\spc

Note that Alice and Bob can both compute the secret commutator: Alice can
premultiply $x^y$ by $x^{-1}$ and Bob can premultiply $y^x$ by
$y^{-1}$ and then compute the inverse:
$[x,y]={(y^{-1}y^x)}^{-1}$. 

The Anshel et al.\ protocol is far from well specified as it stands. In
particular, we have said nothing about our choice of platform group
$G$. Like Ko et al., Anshel et al.\ proposed using braid groups because of the
existence of efficient normal forms for group elements and because the
conjugacy search problem seems hard. See Myasnikov et
al.~\cite[Chapter~5]{Myasnikov} for a discussion of some of the
properties a platform group should have; they discuss the possibilities
of using the following groups as platform groups: Thompson's
group~$F$, matrix groups, small cancellation groups, solvable groups,
Artin groups and Grigorchuck's group.

\subsection{Replacing conjugation}

The Ko et al.\ scheme used conjugation in place of exponentiation in
the Diffie--Hellman protocol, but there are many other
alternatives. For example, we could define $g^a=\phi(a) g a $ and
$g^b=\phi'(b) g b$ for any fixed functions $\phi:A\rightarrow A$ and
$\phi':B\rightarrow B$ (including the identity maps) and the scheme
would work just as well. More generally, we may replace $a$ and
$\phi(a)$ by unrelated elements from $A$: there are protocols based on
the difficulty of the
\emph{decomposition problem}, namely the problem of finding
$a_1,a_2\in A$ such that $h=a_1ga_2$ where $g$ and $h$ are known. See
Myasnikov et al.~\cite[Chapter~4]{Myasnikov} for a discussion of these
and similar protocols; one proposal we find especially interesting is
the Algebraic Eraser~\cite{AAGL,Kalka}. As an example of such a
protocol, we briefly describe a scheme due to Stickel.

\spc\textbf{The Stickel Key Agreement Protocol~\cite{stickel}.}
Let $G=\mathrm{GL}(n,\F_q)$, and let $g\in G$. Let $a,b$ be elements of $G$
of order $n_a$ and $n_b$ respectively, and suppose that $ab \ne ba$. The
group $G$ and the elements $a,b$ are publicly known. If Alice and
Bob wish to create a shared key, they can proceed as follows:
\begin{enumerate}
\item Alice chooses integers $l$, $m$ uniformly at random, where $0 < l < n_a$ and $0 < m < n_b$. She sends $u=a^l g b^m$ to Bob.
\item Bob chooses integers $r,s$ uniformly at random, where $0 < r < n_a$ and $ 0 < s < n_b$. He sends $v=a^r gb^s$ to Alice.
\item Alice computes $k_a=a^l v b^m=a^{l+r}gb^{m+s}$. Bob computes
$k_b=a^rub^s=a^{l+r}g b^{m+s}$.
\item The shared key is thus $k= k_a = k_b$.
\end{enumerate}

\subsection{Logarithmic signatures}

There is an alternative approach to generalising the Diffie--Hellman
scheme: to find a more direct generalisation of the DLP for
groups that are not necessarily abelian.

Let $G$ be a finite group, $S \subseteq G$ a subset of $G$ and $s$ a
positive integer.  For all $1 \leq i \leq s$,
let~$A_i=[\alpha_{i1},\ldots ,\alpha_{ir_i}]$ be a finite sequence of
elements of $G$ of length $r_i>1$, and let $\alpha=[A_1,\ldots ,A_s]$
be the ordered sequence of $A_i$.  We say that $\alpha$ is a
\emph{cover} for $S$ if any $h\in S$ can be written as a product $ h=
h_1 \cdots h_s, $ where $h_i = \alpha_{ik_i} \in A_i$.  If such a
decomposition is unique for every $g\in S$, then $\alpha$ is said to
be a \emph{logarithmic signature} for $S$. One natural way to construct
a logarithmic signature for a group $G$ is to take a subgroup chain
$$1=G_0< G_1< \cdots < G_n=G,$$ and let $A_i$ be a complete set of
coset representatives for $G_{i-1}$ in $G_i$.  Then
$\alpha=[A_1,\ldots ,A_n]$ is a logarithmic signature (a so called
\emph{transversal logarithmic signature}) for $G$.

Given an element $h \in S$ and a cover $\alpha$ of $S$, obtaining a
factorisation
\begin{equation}
\label{eqn:factorisation}
h=
\alpha_{1k_1} \cdots \alpha_{sk_s}
\end{equation}
associated with $\alpha$ could well be a hard problem in
general. Indeed, in some situations the problem is a Discrete
Logarithm Problem. For example, let $G$ be generated by an element $g$
of large order, and define $A_{i+1}=[1,g^{2^i}]$. Let $S=\{g^a\mid0\leq
a\leq 2^s\}$. Then the $i$th bit of the discrete logarithm of $h\in S$ is
equal to $1$ if and only if $k_i=2$ in the
factorisation~\eqref{eqn:factorisation}.

Though there are connections with the DLP, logarithmic signatures
cannot be directly used in discrete logarithm based protocols, as
there is no analogue of exponentiation. They were first used by
Magliveras~\cite{Mag} to construct a symmetric cipher known as
Permutation Group Mappings (PGM). The ideas behind PGM have inspired
several public key cryptosystems based on logarithmic signatures. Qu
and Vanstone~\cite{QuVanstone} proposed a scheme (Finite Group
Mappings, or FGM) based on transversal logarithmic signatures in
elementary abelian $2$-groups. Magliveras, Stinson and van
Trung~\cite{MagStinson} developed two interesting schemes based on
finite permutation groups,
\textit{MST}$_1$ and \textit{MST}$_2$.
More recently, a public key cryptosystem based on Suzuki $2$-groups
(known as
\textit{MST}$_3$) has been proposed by Lempken et al.~\cite{Lempken}.

\subsection{Symmetric schemes}

Group theory has mainly been used in proposals of public key
cryptosystems and key exchange schemes, but has also been used in
symmetric cryptography. We have already mentioned the block cipher
PGM~\cite{Mag}. This cipher satisfies some nice algebraic and
statistical properties (such as robustness, scalability and a large key
space; see~\cite{MagMemon}). However, fast implementation becomes an
issue, making it a rather inefficient cipher compared with more
traditional block ciphers. (An attempt was made to improve PGM by
letting the platform group be a 2-group, but again speed remains an
issue~\cite{Canda}.) This subsection contains two more examples of
group theory being used in symmetric cryptography.

A block cipher such as DES~\cite{nist-fips-46} or
AES~\cite{nist-fips-197} can be regarded as a set $\mathcal{S}$ of
permutations on the set of all possible blocks, indexed by the
key. The question as to whether $\mathcal{S}$ is in fact a
group has an impact on the cipher's security in some situations: if
the set was a group, then encrypting a message twice over using the
cipher with different keys would be no more secure than a single
encryption.  Other properties of the group generated by $\mathcal{S}$
are also of interest cryptographically~\cite{hornauer-euro93} and
attacks have been proposed against ciphers that do not satisfy some of
these properties~\cite{kaliski-rivest-sherman-joc,paterson-fse99}
(though good group theoretic properties are not sufficient to
guarantee a strong cipher~\cite{murphy-paterson-wild-joc}). We note
however that computing the group generated by a block cipher is often
very difficult. For instance, it is known that the group generated by
the DES block cipher is a subgroup of the alternating group
$A_{2^{64}}$~\cite{wernsdorf-euro92}, with order greater than $2^{56}$
(and thus $\mathcal{S}$ for DES is not a
group~\cite{campbell-wiener-crypto92,coppersmith-ibm92}); however
little more is known about its structure. 

Block ciphers themselves are often built as iterated constructions of
simpler key-dependent permutations known as \emph{round functions},
and one can study properties of the permutation groups generated by
these round functions. It has been shown, for instance, that the round
functions of both DES and AES block ciphers are even permutations;
furthermore it can be shown that these generate the alternating group
$A_{2^{64}}$ and $A_{2^{128}}$,
respectively. See \cite{Caranti09,Caranti_preprint,sparr-wernsdorf,wernsdorf-euro92,wernsdorf-fse02}.

Hash function design is a second area of symmetric cryptography where
groups have been used in an interesting
way. Recall~\cite[Chapter~7]{Stinson} that a hash function $H$ is a
function from the set of finite binary strings to a fixed finite set
$X$. It should be easy to compute $H(x)$ for any fixed string $x$, but
it should be computationally infeasible to find two strings $x$ and
$x'$ such that $H(x)=H(x')$. Hash functions are a vital component of
many cryptographic protocols, but their design is still not well
understood. The most widely used example of a hash function is \mbox{SHA-1}
(where SHA stands for Secure Hash
Algorithm). See~\cite{nist-fips-180-1} for a description of this hash
function. Security flaws have been found in SHA-1~\cite{ECRYPT};
the more recent SHA-2 family of hash functions~\cite{nist-fips-180-2}
are now recommended. Z\'emor~\cite{Zemor94} proposed using walks
through Cayley graphs as a basis for hash functions; the most well-known 
concrete proposal from this idea is a hash function of Tillich
and Z\'emor~\cite{TillichZemor94}. We think this hash function
deserves further study, despite a recent (and very beautiful)
cryptanalysis due to Grassl et al.~\cite{Grassl}: see Steinwandt et
al.~\cite{Steinwandt00} and the references there for comments on the
security of this hash function, and see Tillich and
Z\'emor~\cite{TillichZemor08} for some more recent literature.

\section{Cryptanalysis}
\label{sec:cryptanalysis}

In this section, we briefly outline some techniques that
have been developed to demonstrate the insecurity of group-based
schemes.

\subsection{Analysis of braid based schemes}

We begin with braid-based schemes. The interested reader is referred
to the comprehensive survey articles by Dehornoy~\cite{Dehornoy} and
Garber~\cite{Garber}.

In 1969, Garside~\cite{Garside} gave the first algorithm to solve the
conjugacy problem in the braid group $B_n$. (The conjugacy problem
asks whether two braids, in other words two elements of the braid
group, are conjugate.) The question of efficiency of Garside's method
lay dormant until the late 1980's. Since then there has been a great
deal of research, significantly motivated by cryptographic
applications, into finding a polynomial time solution to the conjugacy
problem. Given two braids $x,y\in B_n$, Garside's idea was to
construct finite subsets (so called \emph{summit sets}) $I_x, I_y$ of
$B_n$ such that $x$ is conjugate to $y$ if and only if $I_x=I_y$. An
efficient solution to the conjugacy problem via this method would
yield an efficient solution to the conjugacy search problem (and hence
render the braid based protocol of Ko et al.\
theoretically insecure). However, for a given braid $x$, Garside's
summit set $I_x$ may be exponentially large. The challenge has thus
been to prove a polynomial bound on the size of a suitable invariant
set associated with any given conjugacy class. Refinements to the
summit set method (such as the
\emph{super summit set}, \emph{ultra summit set,} and \emph{reduced super
summit set} methods) have been made over the years, but a polynomial
bound remains elusive. Recent focus has been on an efficient solution
to each of the three types of braids: periodic, reducible or
pseudo-Anasov (according to the Nielsen--Thurston classification);
see~\cite{BirmanI,BirmanII,BirmanIII}.

For the purposes of cryptography however, one need not efficiently
solve the conjugacy problem in order to break a braid-based
cryptosystem: one is free to use the specifics of the protocol being
employed; any algorithm only needs to work for a significant
proportion of cases; heuristic algorithms are quite
acceptable. Indeed, Hofheinz and Steinwandt~\cite{Hofheinz} used a
heuristic algorithm to solve the conjugacy search problem with very
high success rates: their attack is based on the observation that
representatives of conjugate braids in the super summit set are likely
to be conjugate by a permutation braid (a particularly simple
braid). Their attack demonstrates an inherent weaknesses of both the
Ko et al.\ protocol and the Anshel et al.\ protocol for random
instances, under suggested parameters. (This has led researchers to
study ways of generating keys more carefully, to try to avoid easy
instances.) Around the same time, several other powerful lines of
attack were discovered, and we now discuss some of the work that has
been done; see Gilman et al.~\cite{Gilman} for another discussion of
these attacks.

\paragraph{Length-based attacks}

Introduced by Hughes and Tannenbaum~\cite{HughesII}, length-based
attacks provide a neat probabilistic way of solving the conjugacy
search problem in certain cases. Suppose we are given an instance of
the conjugacy search problem in $B_n$. So we are given braids
$x$ and $y^{-1}xy$, and we want to find $y$. Let $l:B_n\rightarrow \mathbb{Z}$
be a suitable length function on $B_n$ (for example, the length of the
normal form of an element). If we can write $y=y'\sigma_i$ for some
$i$, where $y'$ has a shorter length than $y$, then
$l(\sigma_iy^{-1}xy\sigma^{-1}_i)$ should be strictly smaller than
$l(\sigma_jy^{-1}xy\sigma^{-1}_j)$ for $j\neq i$. So $i$ can be
guessed, and the attack repeated for a smaller instance $y'$ of
$y$. The success rate of this probabilistic attack depends on the
specific length function employed. For braid groups, there are a
number of suitable length functions that allow this attack to be
mounted. We comment that length-based attacks need to be modified in practice,
to ensure (for example) that we do not get stuck in short loops; see
Garber et al.~\cite{GarberKaplan} and Ruinskiy et
al.~\cite{RuinskiyShamir}. Garber et al.~\cite{GarberKaplan} and
Myasnikov and Ushakov~\cite{MyUsh} contain convincing attacks on both the
Ko et al.\ and Anshel et al.\ protocols using a length-based approach.

\paragraph{Linear algebra attacks}
The idea behind this attack is quite simple: take a linear
representation of the braid group and solve the conjugacy search
problem using linear algebra in a matrix group. There are two well-known 
representations of the braid group: the Burau representation
(unfaithful for $n\geq 5$) and the faithful Lawrence-Krammer
representation. Hughes~\cite{Hughes} and Lee and Lee~\cite{Lee}
provide convincing attacks on the Anshel et al.\ protocol using the
Burau representation, and Cheon and Jun~\cite{Cheon} provide a
polynomial time algorithm to break the Ko et al.\ protocol using the
Lawrence--Krammer representation. Budney~\cite{Budney} studies the
relationship between conjugacy of elements in the braid group and
conjugacy of their images in the unitary group under the
Lawrence--Krammer representation.

\paragraph{Other directions}
There have been many suggestions made to improve the security of
schemes based on the above protocols. Themes range from changing the
underlying problem (and instead investigating problems such as the
decomposition problem, the braid root problem, the shifted conjugacy
problem and more) to changing the platform group (Thompson's group,
polycyclic groups and others have been suggested). Furthermore,
cryptographers have created other cryptographic primitives based on
the conjugacy search problem, for example authentication schemes and
signature schemes. However, there are no known cryptographic
primitives based on any of these ideas that convincingly survive the
above sketched attacks. It seems to be the pattern that `random' or
`generic' instances of either protocol lead to particularly simplified
attacks. See the book by Myasnikov et al.~\cite{Myasnikov} for more on
this.

\subsection{Stickel's scheme}

Stickel's scheme was successfully cryptanalysed by
Shpilrain~\cite{Shpilrain-attack-stickel}. We include a brief
description of this attack as it is particularly simple, and
illustrates what can go wrong if care is not taken in protocol
design. The attack works as follows.  First note that an adversary
need not recover any of the private exponents $l,m,r,s$ in order to
derive the key $k$.  Instead, it suffices upon intercepting the
transmitted messages $u$ and $v$, to find $n\times n$ matrices $x,y\in
G$ such that $$xa=ax,\;yb=by,\;u=xgy.$$ One can then compute
$$xvy=xa^rgb^sy=a^rxgyb^s=a^rub^s=k.$$

It remains to solve these equations for $x$ and $y$.  The equations
$xa=ax$ and $yb=by$ are linear, since $a$ and $b$ are known.  The
equation $u=xgy$ is not linear, but since $x$ is invertible we can
rearrange: $x^{-1}u=gy$, with $g$ and $u$ known.  Since $xa=ax$ if and only if
$x^{-1}a=ax^{-1}$, we write $x_1=x^{-1}$ and instead solve the following
matrix equations involving $x_1$ and $y$:
$$x_1a=ax_1,\; yb=by,\; x_1u=gy.$$ Setting $x_1=gyu^{-1}$ we can eliminate $x_1$ to solve
$$gyu^{-1}a=agyu^{-1},\; yb=by.$$ Now only $y$ is unknown and we have
$2n^2$ linear equations in $n^2$ variables: a heavily overdetermined
system of linear equations, and an invertible matrix $y$
will be easily found.  Shpilrain's attack is specific to the platform
group $GL(n,\mathbb{F}_q)$.  In particular, it uses the fact that $x$
and $u$ are invertible.  Thus to thwart this attack, it makes sense to
restrict the protocol to non-invertible matrices (since there is no
inversion operation in the key setup).  However, it is unclear whether
or not this actually enhances the security of the protocol.

\subsection{Analysis of schemes based on logarithmic signatures}

How can secure logarithmic signatures be generated? The main problem
with the overwhelming majority of schemes based on logarithmic signatures
is a failure to specify how this should be done. (The Qu--Vanstone
scheme~\cite{QuVanstone} is well specified, but Blackburn, Murphy and
Stern~\cite{BlackburnMurphy} showed this scheme is insecure.)
Magliveras et al.~\cite{MagStinson} had the idea of restricting
the logarithmic signature used in \textit{MST}$_1$ to be
\emph{totally non-transversal}, that is a logarithmic signature $\alpha$ for a group $G$
in which no block $A_i$ of $\alpha$ is a coset of a non-trivial
subgroup of $G$. However, this condition was shown to be insufficient
by Bohli et al.~\cite{Bohli}, who constructed instances of totally
non-transversal logarithmic signatures that were insecure when used
in \textit{MST}$_1$. Key generation is also a problem for
$\mathit{MST}_2$; see~\cite{GonzalezSteinwandt} for a critique of
this.  As for $\mathit{MST}_3$, this was recently cryptanalysed by the
authors~\cite{BlackburnCid}.  Thus it seems that a significant new
idea in this area is needed to construct a secure public key
cryptosystem from logarithmic signatures.

\section{Next Steps}
\label{sec:nextsteps}

Despite ten years of strong interest in group-based cryptography, a
well-studied candidate for a secure, well-specified and efficient
cryptosystem is yet to emerge: schemes that are more `number theoretic' (such
as those based on the elliptic curve DLP) currently have so many
advantages. This is a disappointment (for the group
theorist). However, we do not want to be overly pessimistic: we hope
that the reader is already convinced that the protocols of Ko et al.\
and of Anshel et al.\ are elegant ideas, just waiting for the right
platform group. \emph{Can such a platform group be found?} We need a
candidate group whose elements can be manipulated and stored
efficiently, and an associated problem that is hard in the
overwhelming majority of instances. There has been a great deal of
attention on infinite groups (such as braid groups) that can be
defined combinatorially, but we feel that finite groups deserve a much
closer study; many difficulties disappear when we use finite groups. Note
that groups with small linear representations are often problematic,
as linear algebra can be used to attack such groups; groups with many normal subgroups (such as
$p$-groups, for example) are often
vulnerable to attacks based on reducing a problem to smaller
quotients; groups with permutation representations of low degree are
vulnerable to attacks based on the well developed theory of
computational permutation group theory. So great care must be taken in
the choice of group, and the choice of supposedly hard problem.  More
generally, we can move beyond the Ko et al.\ and Anshel at al.\
schemes, and ask: \emph{Is there a secure and efficient key exchange
protocol based on group theoretic ideas?} There are regular proposals, but the field is still waiting for a proposal that stands up to long-term scrutiny.

We would like to point out that group-based cryptography motivates
some beautiful and natural questions for the pure group theorist. Most
obviously, the cryptosystems above motivate problems in computational
group theory, especially combinatorial group theory. But we would like
to highlight two more problems as examples of the kind of questions
that can arise.

\paragraph{Generic properties}
The cryptosystems described in this survey require that elements and
subgroups of a group $G$ are generated \emph{at random}. This needs to
be defined precisely for this to make sense; one common method would
be to select at random a sequence of integers $\{a_1, a_2, \ldots ,
a_l \}$ of length $l$, and for each $1 \leq i \leq l$, select at random
a generator $x_i$ of $G$. We then output the
\emph{random} element $w = x_1^{a_1} x_2^{a_2} \cdots x_l^{a_l}$.
Many cryptosystems run into problems because randomly generated sets
of elements in the platform group behave in a straightforward way when
$l$ is large. This motivates the study of \emph{generic} properties of
groups, namely properties that hold with probability tending to $1$ as
$l\rightarrow\infty$. For example, Myasnikov and
Ushakov~\cite{MyUsh-JMC} have shown that pure braid groups $PB_n$ have
the strong generic free group property: for any generating set of
$PB_n$, when any $k$ elements are chosen at random as above they
freely generate a free group of rank $k$ generically. An interesting
and natural open problem is: does the same property hold for the braid
groups $B_n$? See Myasnikov et al.~\cite{Myasnikov} for a discussion
of this and related issues.

\paragraph{Short logarithmic signatures}
Let $G$ be a finite group of order $\prod_{j=1}^tp^{a_j}_j$, with
$p_j$ distinct primes.  Let $\alpha=[A_1,\ldots ,A_s]$ be a
logarithmic signature for $G$, with $|A_i|=r_i$ for $1\leq i\leq s$.
Define the \emph{length} of $\alpha$ to be
$l(\alpha):=\sum_{i=1}^{s}r_i$. The length of $\alpha$ is an efficiency measure: it is the number of elements that must be stored in order to specify a typical logarithmic signature of this kind.  Since $|G|=\prod_{i=1}^sr_i$, we must have that $l(\alpha)\geq \sum_{j=1}^ta_jp_j$.  A
logarithmic signature achieving this bound is called a
\emph{minimal logarithmic signature} for $G$. An attractive open problem is: 
does every finite group have a minimal logarithmic signature?  Now, if
$G$ has a normal subgroup $N$ with $G/N\cong H$ and $H$ and $N$ both
have minimal logarithmic signatures then $G$ has a minimal logarithmic
signature. In particular, it is clear that any soluble group has a
minimal logarithmic signature. Moreover, to answer the question in the
affirmative it suffices to consider simple groups only. Minimal
logarithmic signatures have been found for $A_n$, $\mathrm{PSL}_n(q)$,
some sporadic groups and most simple groups of order up to
${10}^{10}$;
see~\cite{Vasco03,GonzalezSteinwandt,Holmes,Lempken05,Magliveras02} for further details.

\spc
Why do we attempt to propose new cryptosystems, when elliptic curve
DLP systems work well? A major motivation is the worry that a good
algorithm could be found for the elliptic curve DLP. This worry has
increased, and the search for alternative cryptosystems has become
more urgent, with the realisation that quantum computers can
efficiently solve both the integer factorisation problem and the
standard variants of the DLP~\cite{Shor}. If quantum computers of a
practical size can be constructed, classical public key cryptography
is in trouble. Cryptosystems, including group-based examples, that are
not necessarily vulnerable to the rise of quantum computers have
become known as~\emph{post-quantum cryptosystems}. A well known
example, invented well before quantum computers were considered, is
the McEliece cryptosystem~\cite{McEliece} based on the difficulty of
decoding error correcting codes. Other examples include lattice-based
cryptosystems (such as the GGH cryptosystem~\cite{GGH,Nguyen}) and
cryptosystems based on large systems of multivariate polynomial
equations (such as the HFE family of
cryptosystems~\cite{KipnisShamir,Patarin}). Though many of these
cryptosystems suffer from having large public keys, they are often
computationally efficient and so we feel that these schemes are more
likely than group-based cryptosystems to produce protocols that will
be used in practice. For a good and recent survey of the area, that
includes more details on all the cryptosystems mentioned above, see
Bernstein et al.~\cite{postquantum}.

We hope the reader is keen to learn more after finishing this
introduction.  We recommend consulting the IACR Cryptology ePrint
Archive~\cite{eprint} or Cornell University's arXiv~\cite{arxiv}
(especially the group theory and cryptography sections) for new
papers; we currently find the ePrint archive the most reliable source
of high quality cryptography. Boaz Tsaban's CGC
Bulletin~\cite{cgc-bulletin} provides regular updates on the main
articles and events in the area. There are many conferences dealing
with cryptographic issues, see~\cite{iacr-events} for a good list;
those conferences sponsored by the IACR are regarded in the field as
being of top quality, though good conferences are not limited to IACR
sponsored events. The \emph{Journal of Cryptology} and \emph{IEEE
Trans.\ Inform.\ Theory} publish excellent papers in the
area; \emph{Designs, Codes and Cryptography} is a well-established source. New
specialist journals that publish papers on group-based cryptography
include the \emph{Journal of Mathematical Cryptology} and
\emph{Groups-Complexity-Cryptology}. For information on group-based schemes based on combinatorial group theory in particular, we would encourage the reader to consult the textbook of Myasnikov et al.~\cite{Myasnikov}.

\paragraph{Acknowledgements} The third author was supported by E.P.S.R.C. PhD studentship EP/P504309/1.

\end{document}